\documentclass[10pt]{amsart2000}

\usepackage{amsmath2000,amssymb,amscd2000,latexsym}

\newcommand{\C}{{\mathbb{C}}}

\newcommand{\N}{{\mathbb{N}}}

\newcommand{\R}{{\mathbb{R}}}

\newcommand{\Z}{{\mathbb{Z}}}

\newcommand{\Bh}{{\mathcal B}}
\newcommand{\Ch}{{\mathcal C}}

\newcommand{\Hh}{{\mathcal H}}

\newcommand{\Kh}{{\mathcal K}}

\newcommand{\Oh}{{\mathcal O}}
\newcommand{\Rh}{{\mathcal R}}

\newcommand{\Th}{{\mathcal T}}

\newcommand{\be}{\mathbf{1}}

\newcommand{\cpr}{\mathrm{cpr }\,}

\newcommand{\dr}{\mathrm{dr}\,}

\newcommand{\halb}{\frac{1}{2}}
\newcommand{\her}{\mathrm{her}}
\newcommand{\hr}{\mathrm{hr}\,}
\newcommand{\id}{\mathrm{id}}
\newcommand{\ord}{\mathrm{ord}\,}

\newcommand{\Prim}{\mathrm{Prim }}

\newcommand{\Sd}{\mathrm{Sd }\, }

\newcommand{\verk}{\mbox{\scriptsize $\,\circ\,$}}

\newcounter{number}[subsection]
\newcounter{altnumber}[section]

\newenvironment{nummer}{\refstepcounter{number}{\noindent\arabic{section}.\arabic{subsection}.\arabic{number}}}{}
\newenvironment{altnummer}{\refstepcounter{altnumber}{\noindent\arabic{section}.\arabic{altnumber}}}{}

\newcommand{\bn}{\noindent\begin{nummer} \rm}
\newcommand{\en}{\end{nummer}}

\newcommand{\altbn}{\noindent \begin{altnummer} \rm}
\newcommand{\alten}{\end{altnummer}}

\newenvironment{ntheorem}{\noindent {\sc Theorem:} \it}{}
\newenvironment{nlemma}{\noindent {\sc Lemma:} \it}{}
\newenvironment{nprop}{\noindent {\sc Proposition:} \it}{}
\newenvironment{ndefn}{\noindent {\sc Definition:} \it}{}
\newenvironment{ncor}{\noindent {\sc Corollary:} \it}{}

\newenvironment{nremarks}{\noindent {\sc Remarks: }}{}
\newenvironment{nexamples}{\noindent {\sc Examples:} }{}

\newenvironment{nproof}{\noindent {\sc Proof:}}{\mbox{}\hfill$\Box$}

\begin{document}

\title{{\sc Covering Dimension and Quasidiagonality}}

\author{Eberhard Kirchberg}
\address{Humboldt--Universit\"at zu Berlin\\ Institut f\"ur Mathematik\\ Unter den Linden 6\\ D-10099 Berlin\\ Germany}

\curraddr{{\sc Mathematisches Institut der Universit\"at M\"unster\\ 
Einsteinstr. 62\\ D-48149 M\"unster\\ Germany}}

\email{kirchbrg@mathematik.hu-berlin.de, kirchbee@math.uni-muenster.de}

\author{Wilhelm Winter}
\address{Mathematisches Institut der Universit\"at M\"unster\\ 
Einsteinstr. 62\\ D-48149 M\"unster\\ Germany}
\email{wwinter@math.uni-muenster.de}

\date{June 2002}


\begin{abstract}
We introduce the decomposition rank, a notion of covering dimension for nuclear $C^*$-algebras. The decomposition rank generalizes ordinary covering dimension and has nice permanence properties; in particular, it behaves well with respect to direct sums, quotients, inductive limits, unitization and quasidiagonal extensions. Moreover, it passes to hereditary subalgebras and is invariant under stabilization. It turns out that the decomposition rank can be finite only for strongly quasidiagonal $C^*$-algebras and that it is closely related to the classification program.
\end{abstract}

\maketitle

\setcounter{section}{-1}

\section{Introduction}

It is well-known that a $C^*$-algebra $A$ is nuclear iff it has the completely positive approximation property, i.e., there is a net $((F_\lambda,\psi_\lambda, \varphi_\lambda))_\Lambda$ with finite-dimensional $C^*$-algebras $F_\lambda$ and completely positive (c.p.) maps $A \stackrel{\psi_\lambda}{\longrightarrow} F_\lambda \stackrel{\varphi_\lambda}{\longrightarrow} A$ such that $\varphi_\lambda \verk \psi_\lambda$ converges to $\id_A$ in the point norm topology. One may always assume the maps $\psi_\lambda$ and $\varphi_\lambda$ to be completely positive contractive (c.p.c.). See \cite{Ta} and \cite{Pa} for general accounts on nuclear $C^*$-algebras and c.p.\ maps.\\
Although all information on $A$ is contained in the approximating system, it is in general not so easy to read off structural properties or $K$-theoretic invariants from $((F_\lambda,\psi_\lambda, \varphi_\lambda))_\Lambda$. One step towards a better understanding of systems of c.p.\ approximations was the introduction of generalized inductive limits and (strong) $NF$ algebras in \cite{BK1}. These concepts are closely related to quasidiagonality (cf.\ \cite{Vo2}), which in turn plays an important r\^ole in Elliott's classification program (see \cite{Ro} for an introduction). The present paper also fits into that framework. Indeed, one motivation for introducing yet another notion of noncommutative covering dimension was, to find conditions under which generalized inductive limit $C^*$-algebras are accessible to methods from $K$-theory, or even to classification.

The triples $(F_\lambda,\psi_\lambda, \varphi_\lambda)$ in some sense are analogues of open coverings of topological spaces. By imposing certain conditions on the maps $\varphi_\lambda$, this point of view was already exploited in \cite{Wi1} and \cite{Wi2} to define the completely positive and the homogeneous rank of a nuclear $C^*$-algebra. Below we modify these conditions to assign to a nuclear $C^*$-algebra $A$ its decomposition rank, $\dr A$. We show that the decomposition rank generalizes topological covering dimension and that it has nice permanence properties. In particular, it behaves nicely with respect to direct sums, quotients, limits, hereditary subalgebras, unitization, stabilization and quasidiagonal extensions. Furthermore it turns out that, if $\dr A$ is finite, then not only $A$, but also every quotient and every representation of $A$ are quasidiagonal.\\
We consider various examples, such as $AF$ algebras (which are precisely the algebras with decomposition rank zero), $AT$ algebras (e.g.\ noncommutative tori), continuous trace algebras, subhomogeneous algebras and certain crossed products.\\
In subsequent work it will be shown that, under a number of extra conditions, simple $C^*$-algebras with finite decomposition rank are classifiable in the sense of \cite{Li}; see Section 6 for a precise statement.

We would like to thank J.\ Cuntz, S.\ Eilers and M.\ R{\o}rdam for several comments and fruitful discussions. The research for this article was supported by the Deutsche Forschungsgemeinschaft.

\section{Covering dimension and graph colourings}

In this section we discuss an alternative characterization of topological covering dimension for normal spaces; this will serve as a commutative model for our decomposition rank.

\altbn
For the convenience of the reader we recall the following characterization of covering dimension (cf.\ \cite{HW}, Theorem V.8; there, the space was assumed to be separable and metrizable):

\begin{ndefn}Let $X$ be a topological space.\\
a) The order of a family $(U_{\lambda})_{\Lambda}$ of subsets of $X$ does not exceed $n$, if for any $n +2 $ distinct indices $\lambda_0 , \ldots , \lambda_{n+1} \in \Lambda$ we have $\bigcap^{n+1}_{i=0} U_{\lambda_i} = \emptyset$.\\
b) The covering dimension of $X$ does not exceed $n$, $\dim X \le n$, if every finite open covering of $X$ has an open refinement of order not exceeding $n$. We say $\dim X = n$, if $n$ is the least integer such that $\dim X \le n$.
\end{ndefn}
\alten

\altbn{\label{complexes}}
A finite simplicial complex $\Sigma$ is a collection of subsets (called faces) of a finite vertex set $V(\Sigma)$ such that, if $\sigma \in \Sigma$, then $\sigma' \subset \sigma$ implies $\sigma' \in \Sigma$. Let $|\Sigma|$ be the geometric realization of $\Sigma$. Choosing a bijection $V(\Sigma) \to \{1, \ldots, k\}$ (where $k := {\mathrm{ card }}\, V(\Sigma)$), $|\Sigma|$ can be obtained as a compact subset of the standard simplex $\Delta^{k-1}$ in $\R^k$. There is a canonical open covering $(A_\gamma)_{\gamma \in V(\Sigma)}$ of $|\Sigma|$ which comes from open stars around vertices in $\Delta^{k-1}$. More precisely, $A_\gamma$ consists of those points of $|\Sigma|$ for which the coordinate belonging to $\gamma$ vanishes.\\
For each $\sigma \in \Sigma$, let $\dim \sigma$ be the combinatorial dimension of $\sigma$; this is ${\mathrm{ card }}\, \sigma -1$. Set $\dim \Sigma := \max \{ \dim \sigma \, | \, \sigma \in \Sigma \}$; this coincides with the covering dimension of $|\Sigma|$.\\
One can form the barycentric subdivision $\Sd \Sigma$ of $\Sigma$ as follows: Let $V(\Sd \Sigma) := \Sigma$ and define a subset $S \subset \Sigma$ to be in $\Sd \Sigma$ iff $S$ can be written as $\{ \sigma_1, \ldots, \sigma_k \}$ such that $\sigma_i \subset \sigma_{i+1} \in \Sigma$, $i=1, \ldots, k-1$. Then $\Sd \Sigma$ is a simplicial complex and an affine homeomorphism $|\Sd \Sigma| \to |\Sigma|$ is given by sending each vertex to the barycenter of the corresponding face in $|\Sigma|$.\\
If $\Gamma$ is a graph with vertex set $V(\Gamma)$ and edge set $E(\Gamma)$, we say $\Gamma$ is $n$-colourable if there is an $n$-colouring, i.e.\ a map $c: V(\Gamma) \to \{0, \ldots, n-1\}$ such that $c(\gamma) \neq c(\gamma')$ whenever $\gamma \neq \gamma'$ and $\{\gamma, \gamma'\} \in E(\Gamma)$.\\
The 1-skeleton of $\Sigma$ is the subcomplex given by nonempty subsets of $V(\Sigma)$ which are in $\Sigma$ and contain no more than two elements. This may be regarded as a graph with vertex set $V(\Sigma)$ in the obvious way.
\alten

\altbn{\label{1-skeleton}}
\begin{nprop}
Let $\Sigma$ be a simplicial complex with $\dim \Sigma = n$. Then the 1-skeleton of $\Sd \Sigma$, regarded as a graph, is $(n+1)$-colourable.
\end{nprop}

\begin{nproof}
Define a map $c : V(\Sd \Sigma) = \Sigma \to \{0, \ldots, n\}$ by $c(\sigma) := \dim \sigma$. If $\{\sigma, \sigma'\}$ is in the 1-skeleton of $\Sd \Sigma$, then $\sigma \subset \sigma'$ or $\sigma' \subset \sigma$. But if $\sigma \neq \sigma'$, this means that their cardinalities are different, so $c(\sigma) = \dim \sigma \neq \dim \sigma' = c(\sigma')$.
\end{nproof}
\alten

\altbn{\label{decomposable-covering}}
\begin{ndefn}
Let $X$ be a set and $(U_\lambda)_\Lambda$ a collection of subsets. $(U_\lambda)_\Lambda$ is $n$-decomposable, if there is a decomposition $\Lambda = \coprod_{j=0}^n \Lambda_j$ such that
\[
(\lambda, \lambda' \in \Lambda_j \Rightarrow U_\lambda \cap U_{\lambda'} = \emptyset) \, \forall \, j \in \{0, \ldots, n \} \, .
\]
\end{ndefn}
\alten

\altbn{\label{decomposable=colourable}}
\begin{nremarks}
(i) It is trivial that, if $(U_\lambda)_\Lambda$ is $n$-decomposable, then it has order not exceeding $n$; even more, it has strict order (in the sense of \cite{Wi1}, Definition 2.7) not exceeding $n$.\\
(ii) From $(U_\lambda)_\Lambda$ as above we can obtain a graph $\Gamma$ by setting $V(\Gamma) := \Lambda$ and $E(\Gamma) := \{ \{ \lambda, \lambda' \} \subset \Lambda \, | \, \lambda \neq \lambda', \, U_\lambda \cap U_{\lambda'} \neq \emptyset \}$. Then $(U_\lambda)_\Lambda$ is $n$-decomposable iff $\Gamma$ is $(n+1)$-colourable. In particular note that, if $\Sigma$ is a simplicial complex with $\dim \Sigma = n$ and $(A_\sigma)_\Sigma$ is the open covering of $|\Sd \Sigma|$ given by open stars around vertices in $\Sd \Sigma$, then $(A_\sigma)_\Sigma$ is $n$-decomposable, since the underlying graph coincides with the 1-skeleton of $\Sd \Sigma$, which is $(n+1)$-colourable by Proposition \ref{1-skeleton}.
\end{nremarks}
\alten

\altbn{\label{dim=decomposable}}
The definition of our decomposition rank is motivated by the following characterization of covering dimension. We note that this is implicitly contained in the proof of \cite{HW}, Theorem V 8 (at least for separable metric spaces), which is a consequence of the decomposition theorem for covering dimension (\cite{HW}, Theorem III.3). For completeness we include an explicit (and different) proof which makes use of barycentric subdivisions.
 
\begin{nprop}
For a normal space $X$ the following are equivalent:\\
(i) $\dim X \le n$\\
(ii) every finite open covering of $X$ has an $n$-decomposable finite open refinement.
\end{nprop}

\begin{nproof}
(ii) $\Rightarrow$ (i) is obvious; we show (i) $\Rightarrow$ (ii), which is a variation of \cite{Wi1}, Proposition 2.8. \\
Given an open covering of $X$, take an open refinement $(V_\lambda)_\Lambda$ of order not exceeding $n$. Choose a partition of unity $(h_\lambda)_\Lambda$ subordinate to this refinement to obtain a map $h : X \to |\Sigma|$ for some $n$-dimensional simplicial complex $\Sigma$ with vertex set $V(\Sigma) = \Lambda$. Let $(A_\sigma)_\Sigma$ be the canonical open covering (given by open stars around vertices, cf.\ \ref{complexes}) of $|\Sd \Sigma|$; by Proposition \ref{1-skeleton} and Remark \ref{decomposable=colourable}(ii) $(A_\sigma)_\Sigma$ is $n$-decomposable.\\
Set $U_\sigma := h^{-1}(A_\sigma) \subset X$, then $(U_\sigma)_\Sigma$ is an open covering of $X$ since $(A_\sigma)_\Sigma$ covers $|\Sd \Sigma|$ and therefore (identifying $|\Sigma|$ and $|\Sd \Sigma|$) also $|\Sigma|$.\\
$(U_\sigma)_\Sigma$ refines the initial covering since it refines $(V_\lambda)_\Lambda$: For each $\sigma \in \Sigma$ choose $\lambda \in \Lambda$ such that $\lambda \in \sigma$, then (again identifying $|\Sd \Sigma|$ and $|\Sigma| \subset \R^\Lambda$) $A_\sigma$ is contained in the open star of $|\Sigma|$ around $\lambda$, which means that the coordinate belonging to $\lambda$ of each point in $A_\sigma$ is nonzero. In other words, we have $h_\lambda(h^{-1}(x)) \neq 0 \; \forall \; x \in A_\sigma$, so $A_\sigma \subset V_\lambda$.\\
Finally, the same decomposition of $\Sigma$ as for $(A_\sigma)_\Sigma$ also shows that $(U_\sigma)_\Sigma$ is $n$-decomposable.
\end{nproof}
\alten

\section{Order zero maps and decomposability}

We will define the decomposition rank in terms of c.p.\ approximations $(F_\lambda, \psi_\lambda, \varphi_\lambda)$ by imposing a certain condition on the maps $\varphi_\lambda$. Below we examine this condition and show that it is stable under slight perturbations. It is this flexibility that will allow us to derive many more permanence properties for the decomposition rank than (at least at the present stage) are known for the homogeneous or the completely positive rank (cf.\ \cite{Wi2}).

\altbn{\label{strict-order}} 
Let $A$ and $F$ be $C^*$-algebras, $F$ finite-dimensional. Recall from \cite{Wi1}, Definition 3.1, that a c.p.\ map $\varphi: F \to A$ is said to have strict order $n$, $\ord \varphi = n$, if $n$ is the least integer such that the following holds:\\
For any set $\{e_0, \ldots, e_{n+1}\} \subset F$ of pairwise orthogonal minimal projections there are $i, \, j \in \{0, \ldots, n+1\}$ such that $\varphi(e_i) \perp \varphi(e_j)$.

Of course strict order zero maps are of particular interest; they can be characterized as follows (cf.\ \cite{Wi1}, Proposition 4.1.1 a)): 

\begin{nprop}
If $\varphi : F \to A$ is c.p.c.\ with $\ord \varphi = 0$, then there is a unique $*$-homomorphism $\pi_\varphi : CF \to A$ such that $\pi_\varphi(h \otimes x) = \varphi(x) \; \forall \, x \in F$, where $CF$ is the cone $\Ch_0((0,1]) \otimes F$ over $F$ and $h := \id_{(0,1]}$ is the generator of $\Ch_0((0,1])$. Conversely, any $*$-homomorphism $\pi : CF \to A$ induces such a c.p.c.\ map $\varphi$ of strict order zero. 
\end{nprop}
\alten

\altbn{\label{d-n-decomposable}}
\begin{ndefn}
Let $A$ be a $C^*$-algebra. A c.p.\ map $\varphi : \bigoplus_{i=1}^s M_{r_i} \to A$ is $n$-decomposable, if there is a decomposition $\{1, \ldots, s\} = \coprod_{j=0}^n I_j$ such that the restriction of $\varphi$ to $\bigoplus_{i \in I_j} M_{r_i}$ has strict order zero for each $j \in \{0, \ldots, n\}$.
\end{ndefn}
\alten

\altbn{\label{l-decomposable}}
\begin{nlemma}
Let $A$, $F$ be $C^*$-algebras, $F= M_{r_1} \oplus \ldots \oplus M_{r_s}$ finite-dimensional, and let $\varphi : F \to A$ be an $n$-decomposable c.p.c.\ map. Then, for any $I \subset \{1, \ldots , s\}$, $\| \sum_{i \in I} \varphi(\be_i) \| \le (n+1) \cdot \max_{i \in I} \{ \| \varphi(\be_i) \| \}$.
\end{nlemma}

\begin{nproof}
Since $\varphi$ is $n$-decomposable, there is a decomposition $I = \coprod_{j=0}^n I_j$ of $I$ such that $\ord (\varphi_{I_j}) = 0 \; \forall \, j$, where $\varphi_{I_j} := \varphi|_{(\bigoplus_{i \in I_j} M_{r_i})}$. But 
\begin{eqnarray*}
\| {\textstyle \sum_{i \in I}} \varphi(\be_i) \| & \le & (n+1) \cdot \max_{j \in \{0, \ldots, n\}} \{ \| {\textstyle \sum_{i \in I_j}} \varphi(\be_i) \| \\
& = & (n+1) \cdot \max_{j \in \{0, \ldots, n\}} \max_{i \in I_j} \{ \| \varphi(\be_i) \| \} \, ,
\end{eqnarray*}
since $\varphi(\be_i) = \varphi_{I_j}(\be_i) \perp \varphi_{I_j}(\be_{i'}) = \varphi(\be_{i'})$ if $i, \, i' \in I_j$.
\end{nproof}
\alten

\altbn{\label{liftable}}
By \cite{Lo}, Theorems 10.1.11 and 10.2.1, cones over finite-dimensional $C^*$-algebras are projective. In connection with Proposition 2.1 this shows that, if $J \lhd A$ is an ideal and $\varphi: F \to A/J$ is c.p.(c.) with $\ord \varphi = 0$, then $\varphi$ has a c.p.(c.) lift $\hat{\varphi} :F \to A$ with $\ord \hat{\varphi} = 0$ (cf.\ \cite{Wi2}, Proposition 1.2.4). If $\varphi$ is $n$-decomposable, we can lift each $\varphi_{I_j}$ (defined as in the proof of Lemma \ref{l-decomposable}) to an order zero map $\hat{\varphi}_{I_j}$ separately. Therefore, $n$-decomposable maps can be lifted to $n$-decomposable maps; however, the lift might be larger in norm.
\alten

\altbn{\label{relations}}
In \cite{Wi2}, 1.2.3, maps of strict order zero were described in terms of finitely many generators and relations. This carries over to $n$-decomposable maps; for completeness, we state the relations explicitly. Let $F = M_{r_1} \oplus \ldots \oplus M_{r_s}$ be fixed. For some decomposition $\{1, \ldots,s\} = \coprod_{j=0}^n I_j$ consider the relations
\[
\begin{array}{lc}
(\Rh) & \left\{
\begin{array}{l}
\|x_1^{(i)}\| \le 1, \; x_1^{(i)} \ge 0, \; \|x_k^{(i)}\| \le 1, \\
\\
x_k^{(i)} x_l^{(i)} = 0 , \; (x_k^{(i)})^* x_l^{(i)}  = \delta_{k,l} \cdot (x_1^{(i)})^2 , \\
\\
{\left( \sum_{m=1}^{r_i} x_m^{(i)} (x_m^{(i)})^* \right) \,  \left( \sum_{m=1}^{(r_{i'})} x_m^{(i')} (x_m^{(i')})^* \right) } = 0 
\end{array}
\right. \\
\\
(\Rh) & \left\| {\sum_{t=1}^{s} \sum_{m=1}^{r_t} x_m^{(t)} (x_m^{(t)})^*  } \right\| \le 1 
\end{array}
\]
for all $j \in \{0, \ldots,n\}, \; i,i' \in I_j$ and $k,l \in \{2, \ldots, r_i \}$. Then a c.p.\ map $\varphi: F \to A$ is $n$-decomposable if and only if there is a set of matrix units $e_{k,l}^{(i)}$, $i \in \{1, \ldots,s\}$ and $k,l \in \{1, \ldots, r_i\}$, and some decomposition $\{1, \ldots,s\} = \coprod_{j=0}^n I_j$ such that the elements $x_1^{(i)} := \varphi(e_{1,1}^{(i)})$ and $x_k^{(i)} := \varphi(e_{k-1,k}^{(i)})$ (for $i \in \{1, \ldots,s\}$ and $k \in \{2, \ldots, r_i \}$) satisfy the relations $(\Rh)$. $\varphi$, additionally, is contractive if and only if the $x_k^{(i)}$ also satisfy $(\Rh1)$.
\alten

\altbn{\label{weakly-stable}}
From \ref{liftable} and \ref{relations} we see that the relations $(\Rh)$ are liftable; by \cite{Lo}, Proposition 10.1.7 and Theorem 4.1.4, this implies that they are also weakly stable. But a small perturbation of the $x_k^{(i)}$ of no more than, say, $\beta$, will disturb $(\Rh1)$ by no more than $\beta \cdot 2 \sum_{i=1}^s r_i$, and this shows that even $(\Rh) \cup (\Rh1)$ are weakly stable, in other words:

\begin{nprop}
For $F = M_{r_1} \oplus \ldots \oplus M_{r_s}$ and $\eta > 0$ there is $\delta > 0$ such that the following holds:\\
Suppose $\varphi : F \to A$ is c.p.\ with $\| \varphi(\be_F) \| < 1 + \delta$ and, for some set 
\[
\{ e^{(i)}_{k,l} \, | \,  i =1, \ldots, s \, ; \, k,l= 1, \ldots, r_i \}
\]
of matrix units for $F$, the elements $x_k^{(i)}$, defined as in \ref{relations}, satisfy $(\Rh)$ and $(\Rh1)$ up to $\delta$. Then there is an $n$-decomposable c.p.c.\ map $\varphi' : F \to A$ with $\| \varphi - \varphi' \| < \eta$. 
\end{nprop}
\alten

\section{Decomposition rank}

We are now prepared to define the decomposition rank and derive its most important permanence properties.

\altbn{\label{d-dr}}
\begin{ndefn}
Let $A$ be a separable $C^*$-algebra. $A$ has decomposition rank $n$, $\dr A = n$, if $n$ is the least integer such that the following holds: Given $\{b_1, \ldots, b_m\} \subset A$ and $\varepsilon > 0$, there is a c.p.\ approximation $(F, \psi, \varphi)$ for $b_1, \ldots, b_m$ within $\varepsilon$ such that $\varphi$ is $n$-decomposable.\\
If no such $n$ exists, we write $\dr A = \infty$.   
\end{ndefn}
\alten

\altbn{\label{r-dr}}
\begin{nremarks}
(i) If, in the preceding definition, $\varphi$ is merely asked to have strict order not exceeding $n$, this defines the completely positive rank, $\cpr A$. If $\ord \varphi \le n$ and $\ord \varphi|_{M_{r_i}}=0 \; \forall \; i$ (with $F= M_{r_1} \oplus \ldots \oplus M_{r_s}$), this defines the homogeneous rank, $\hr A$; see \cite{Wi1} and \cite{Wi2} for details.\\ 
We see that, for any $A$, $\cpr A \le \hr A$ and it is conceivable that we always have equality. In \cite{Wi2} it was proved that $\cpr A = \hr A$ if $A$ is simple. More generally, one can show that, if $A$ has no irreducible representation of dimension less than or equal to $\cpr A + 1$, then $\cpr A = \hr A$.\\
It is easy to see that $n$-decomposable maps have strict order not exceeding $n$. By definition they have strict order zero on  matrix algebras, from which follows that $\hr A \le \dr A$ for any $A$. At the present stage it is not clear wether the completely positive, the homogeneous and the decomposition rank can take different values.\\
(ii) In the definition it suffices to find c.p.\ approximations for normalized positive elements $b_1, \ldots, b_m$. Even more, by perturbing the $b_j$ slightly (using an approximate unit for $A$), one can show that it suffices to approximate subsets $\{h, b_1, \ldots, b_m\}$ of positive normalized elements with the additional property that $h b_j = b_j \; \forall \, j$ (this will be useful if $A$ is not unital).\\
(iii) It is clear that $\dr A \le n$ iff $A$ has a system $(F_\lambda, \psi_\lambda, \varphi_\lambda)_\Lambda$ of c.p.\ approximations such that $\varphi_\lambda$ is $n$-decomposable $\forall \, \lambda$.
\end{nremarks}
\alten

\altbn{\label{dr-permanence}}
Decomposition rank has nice permanence properties. Essentially the same proofs as in \cite{Wi1}, Section 3, show that
\begin{itemize}
\item[(i)] $\dr (A \oplus B) = \max \{\dr A, \dr B\}$
\item[(ii)] $\dr A \le \underline{\lim} \, \dr A_n$ if $A = \lim_\to A_n$
\item[(iii)] $\dr (A/J) \le \dr A$ if $J \lhd A$ is an ideal;
\end{itemize}
furthermore we will see that 
\begin{itemize}
\item[(iv)] $\dr B \le \dr A$ if $B \subset_{\her} A$ is a hereditary subalgebra.
\end{itemize}
The proof of \cite{Wi2}, Proposition 3.1.4 also shows that, for any two $C^*$-algebras $A$ and $B$, 
\[
\dr (A \otimes B) \le (\dr A + 1) \cdot (\dr B +1) -1 \, ;
\]
if $B$ is $AF$ we have $\dr (A \otimes B) \le \dr A$.
\alten

\altbn{\label{dim=dr}}
\begin{nprop}
Let $X$ be a second countable locally compact space. Then $\dr (\Ch_0(X)) = \dim X$.
\end{nprop}

\begin{nproof}
First note that, for $\varphi : \C^k \to \Ch_0(X)$ c.p.c., $\varphi$ is $n$-decomposable as a map iff $(\varphi(e_i)^{-1}((0,\infty)))_{i \in \{1, \ldots, k\}}$ is as a collection of subsets. If $\dim X \le n$, from Proposition \ref{dim=decomposable} we see that each finite open covering has an $n$-decomposable finite open refinement. Now the proof of \cite{Wi1}, Proposition 3.5 (using decomposability instead of strict order) shows that $\dr (\Ch_0(X)) \le \dim X$.\\
We have $\dim X = \cpr (\Ch_0(X)) \le \hr (\Ch_0(X)) \le \dr (\Ch_0(X))$ by Remark \ref{r-dr}(i) and \cite{Wi1}, Proposition 3.18.
\end{nproof}
\alten

\altbn{\label{multiplicative-domain}}
For frequent use we note a straightforward and well-known consequence of Stinespring's theorem (cf.\ \cite{KiR}, Lemma 7.11):

\begin{nlemma}
Let $A$, $B$ be $C^*$-algebras and $\varphi: B \to A$ a c.p.c.\ map. Then, for any $x,y \in B$, we have $\| \varphi(xy) - \varphi(x) \varphi(y)\| \le \| \varphi(xx^*) - \varphi(x) \varphi(x^*)\|^\halb \|y\|$. 
\end{nlemma}
\alten

\altbn{\label{almost-multiplicative}}
\begin{nlemma}
Let $A$ and $B$ be $C^*$-algebras, $a \in A_+$ with $\|a\| \le 1$ and $\eta > 0$. If $A \stackrel{\psi}{\longrightarrow} B \stackrel{\varphi}{\longrightarrow} A$ are c.p.c.\ and $\| \varphi \psi (a) - a \| , \, \| \varphi \psi(a^2) - a^2 \| < \eta$, then, for all $b \in B$, $\| \varphi (\psi(a)b) - \varphi \psi (a) \varphi (b) \| < 3^\halb \eta^\halb {\|b\|}$.
\end{nlemma}

\begin{nproof}
We have $\| \varphi (\psi(a)^2) - (\varphi \psi(a))^2 \| \le \| \varphi \psi(a^2) - \varphi \psi(a)^2 \| \le 3 \eta$, so $\| \varphi (\psi(a) b) - \varphi \psi(a) \varphi(b) \| < {(3 \eta)^\halb}{\|b\|} \; \forall \, b \in B$ by Lemma \ref{multiplicative-domain}.
\end{nproof}
\alten

\altbn{\label{functions}}
Recall the following notation from \cite{Wi2}, 3.2.4: For positive numbers $\alpha$ and $\varepsilon$ define continuous positive functions on $\R$
\[
f_{\alpha,\varepsilon} (t) :=  \left\{ 
    \begin{array}{cl}
0 & \mbox{for} \; t \le \alpha \\
t & \mbox{for} \; \alpha + \varepsilon \le t \\
\mbox{linear} &  \mbox{elsewhere}
\end{array} \right.
\]
and
\[
g_{\alpha , \varepsilon} (t) :=  \left\{ 
    \begin{array}{cl}
0 & \mbox{for} \; t \le \alpha \\
1 & \mbox{for} \; \alpha + \varepsilon \le t \\
\mbox{linear} & \mbox{elsewhere \,;}
\end{array} \right.
\]
let $\chi_{\alpha}$ denote the characteristic function of $[\alpha , \infty)$.
\alten

\altbn{\label{hereditary-perturbation}}
\begin{nlemma}
Let $A$, $F$ be $C^*$-algebras, $F$ finite-dimensional, and let $\varphi : F \to A$ be c.p.c.\ with $\ord \varphi = 0$. Suppose that there is $h \in A_+$, $\|h\| \le 1$, satisfying $\|\varphi(\be_F) - \varphi(\be_F) h\| < \delta$ for some $0 < \delta < \frac{1}{16}$.\\
Then there is a c.p.c.\ map $\hat{\varphi} : F \to \overline{hAh}$ with $\ord \hat{\varphi} = 0$ and such that $\| \hat{\varphi} (x) - \varphi(x) \| < 8 \delta^{\frac{1}{2}} \; \forall \, x \in F_+$ with $\|x\| \le 1$.
\end{nlemma}

\begin{nproof}
Let $b:= \varphi(\be_F)$, then by (the proof of) \cite{Wi1}, Proposition 4.4.1 a), there is a $*$-homomorphism $\rho : F \to A''$ such that $b \in \rho(F)'$ and $b \rho(x) = \varphi(x) \; \forall \, x \in F$. Set $b':= f_{\delta^\halb, 2 \delta^\halb} (b)$ and $d:= g_{0, \delta^\halb}(b)$; we have $d b' = b'$, $\|b' - b\| \le \delta^\halb$ and
\[
\| d (\be - h)\|^2 = \| (\be -h)d^2(\be -h) \| \le \delta^{-1} \|(\be -h)b^2(\be -h)\| \le \delta \, ,\]
so $\|d (\be -h)\| \le \delta^\halb$. Writing $b'' := (b')^\halb h^2 (b')^\halb$ we obtain
\begin{eqnarray*}
b' & \ge & b''\\
& = & (b')^\halb d h^2 d (b')^\halb\\
& \ge & (b')^\halb (d^2 - 2 \delta^\halb \cdot \be) (b')^\halb\\
& = & (1 - 2 \delta^\halb) \cdot b' \, .
\end{eqnarray*}
This means we have $\|b' - b''\| \le 2 \delta^\halb$ and in the same way it follows that $\| (b')^\halb - (b'')^\halb \| \le 2 \delta^\halb$. Furthermore we see that $b'$ and $b''$ have the same support projections in $A''$, which in turn means that they generate the same hereditary subalgebras. In particular we obtain
\[
C:= \overline{(b')^\halb h A h (b')^\halb} = \overline{b'' A b''} = \overline{b' A b'} \subset A
\]
and
\[
\tilde{C} := \overline{(b')^\halb h A'' h (b')^\halb} = \overline{b'' A'' b''} = \overline{b' A'' b'} \subset A'' \, .
\]
Let $u \in A''$ be the partial isometry in the polar decomposition 
\[
h (b')^\halb = u ((b')^\halb h^2 (b')^\halb)^\halb = u (b'')^\halb \, ,
\]
then $u$ induces an isomorphism 
\[
\pi : \tilde{C} \stackrel{\cong}{\longrightarrow} \overline{h (b')^\halb A'' (b')^\halb h} \, , \; y \mapsto u y u^* \, .
\]
Note that $\pi$ restricts to an isomorphism $C \cong \overline{h (b')^\halb A (b')^\halb h} \subset \overline{h A h}$ and that 
\def\theequation{$\ast$}
\begin{equation}
\pi((b'')^\halb y (b'')^\halb) = h (b')^\halb y (b')^\halb h \; \forall \, y \in A'' \, .
\end{equation}
Define $\hat{\varphi} : F \to \overline{h A h}$ by $\hat{\varphi}(x) := \pi(b' \rho(x))$ for $x \in F$; $\hat{\varphi}$ is well-defined, since $b' \rho(x) = (b')^\halb \rho(x) (b')^\halb \in C$, it clearly is c.p.\ and contractive. $\hat{\varphi}$ has strict order zero since $b' \rho ( \, . \, )$ has and $\pi$ is a $*$-homomorphism. Finally, for $x \in F_+$ with $\|x\| \le 1$ we have
\begin{eqnarray*}
\| \hat{\varphi}(x) - \varphi(x)\| & = & \| \pi((b')^\halb \rho(x) (b')^\halb) - \varphi(x) \| \\
& \le & \| \pi((b')^\halb \rho(x) (b')^\halb) - \pi((b'')^\halb \rho(x) (b'')^\halb) \| \\
& & + \| \pi((b'')^\halb \rho(x) (b'')^\halb) - (b')^\halb \rho(x) (b')^\halb \| \\
& & + \| (b')^\halb \rho(x) (b')^\halb -b^\halb \rho(x) b^\halb \| \\
& < & 4 \delta^\halb + 2 \delta^\halb + 2 \delta^\halb \\
& = & 8 \delta^\halb \, ,
\end{eqnarray*}
where we used ($\ast$) and
\[
\| (b')^\halb (\be - h)\| = \| (b')^\halb d (\be-h)\| < \delta^\halb
\]
for the second estimate.
\end{nproof}
\alten

\altbn{\label{dominated-approximations}}
\begin{nlemma}
Let $A$ be a separable $C^*$-algebra with $\dr A=n < \infty$. Suppose $0 < \eta \le 1$ is given and $h, b_1, \ldots, b_m \in A_+$ are normalized elements that satisfy $h b_j = b_j \; \forall \, j$.\\
Then there are a c.p.\ approximation $(F, \psi, \varphi)$ for $\{h, b_1, \ldots, b_m\}$ within $\eta$ and projections $p^{(l)} \le q^{(l)} \in F$, $l=0, \ldots, n$, with the following properties:
\begin{itemize}
\item[(1)] the $q^{(l)}$ are pairwise orthogonal, central in $F$ and sum up to $\be_F$
\item[(2)] $\ord (\varphi|_{q^{(l)}F}) = 0 \; \forall \, l$ (so $\varphi$ is $n$-decomposable)
\item[(3)] for all $j,l$ and $p:= \sum_{l=0}^n p^{(l)}$ the numbers 
\[
\begin{array}{l}
\| \varphi(p \psi(b_j) p) - b_j\| \, , \\
\| \varphi (q^{(l)} \psi(b_j) q^{(l)} - p^{(l)} \psi(b_j) p^{(l)})\|\, , \\
\| \varphi(p^{(l)}) - \varphi(p^{(l)}) h\| \, , \\
\| \varphi(p^{(l)}) - \varphi(p^{(l)}) \varphi(\be_F) \| \\
\mbox{and } \\
\| \varphi (p^{(l)} - p^{(l)} \psi(h)) \|
\end{array}
\]
are all smaller than $\eta$.
\end{itemize}
\end{nlemma}

\begin{nproof}
Choose a system $(F_\nu, \psi_\nu, \varphi_\nu)_\N$ of c.p.\ approximations for $A$ with $\varphi_\nu$ $n$-decomposable for each $\nu$. Choosing a free ultrafilter $\omega$ on $\N$, we may consider the ultrapowers (cf.\ \cite{Ro}, 6.2) and the induced maps 
\[
\textstyle
\prod_\omega A \stackrel{\psi_\omega}{\longrightarrow} \prod_\omega F_\nu \stackrel{\varphi_\omega}{\longrightarrow} \prod_\omega A \; ,
\]
then  $\varphi_\omega \psi_\omega \iota (a) = \iota(a) \in A_\omega \; \forall \; a \in A$, where $\iota : A \hookrightarrow A_\omega$ is the natural embedding (we will omit the $\iota$ in what follows). Note that 
\[
\varphi_\omega \psi_\omega (a^*) \varphi_\omega \psi_\omega (a) \le \varphi_\omega (\psi_\omega (a^*) \psi_\omega(a)) \le \varphi_\omega \psi_\omega (a^* a) \; \forall \, a \, ,
\]
so $\varphi_\omega$ is multiplicative on $C^*(\psi(A))$. But then Lemma \ref{multiplicative-domain} yields 
\[
\varphi_\omega(x \psi_\omega(a)) = \varphi_\omega(x) \varphi_\omega \psi_\omega (a) = \varphi_\omega(x) a \;  \forall \, x \in {\textstyle \prod_\omega F_\nu}, \, a \in A \, .
\]
For each $\nu \in \N$ define $p_\nu := \chi_{1-\frac{\eta^2}{4}} (\psi_\nu(h)) \in F_\nu$ and denote the image of $(p_\nu)_\N$ in $\prod_\omega F_\nu$ by $p_\omega$;  we clearly have $p_\omega g_{1-\frac{\eta^2}{4},1} (\psi_\nu(h)) = g_{1-\frac{\eta^2}{4},1} (\psi_\nu(h))$ (cf.\ \ref{functions} for the definition of $\chi_\alpha$ and $g_{\alpha, \beta}$). Now for each $j$ we obtain
\begin{eqnarray*}
\varphi_\omega(p_\omega) b_j & = & \varphi_\omega(p_\omega) g_{1-\frac{\eta^2}{4},1} (h) b_j\\
& = & \varphi_\omega(p_\omega) g_{1-\frac{\eta^2}{4},1} (\varphi_\omega \psi_\omega(h)) b_j\\
& = & \varphi_\omega(p_\omega) \varphi_\omega (g_{1-\frac{\eta^2}{4},1} (\psi_\omega(h))) b_j\\
& = & \varphi_\omega(p_\omega g_{1-\frac{\eta^2}{4},1} (\psi_\omega(h))) b_j\\
& = & \varphi_\omega(g_{1-\frac{\eta^2}{4},1} (\psi_\omega(h))) b_j\\
& = & g_{1-\frac{\eta^2}{4},1} (h) b_j \\
& = & b_j \, .
\end{eqnarray*}
Furthermore we have $\varphi_\omega (p_\omega \psi_\omega (h)) = \varphi_\omega(p_\omega)h$ and $\| p_\omega \psi_\omega(h) - p_\omega \| \le \frac{\eta^2}{4}$, hence 
\[
\| \varphi_\omega(p_\omega) - \varphi_\omega(p_\omega)h \| \le \frac{\eta^2}{4} \, . 
\]
It is therefore possible to choose $\nu \in \N$ such that, for $j = 1,\ldots,m$, we have
\def\theequation{$\ast \ast$}
\begin{equation}
\left.
\begin{array}{l}
\| \varphi_\nu \psi_\nu (h) - h\|  \\
\| \varphi_\nu \psi_\nu (b_j) - b_j \| \\
\| \varphi_\nu( p_\nu \psi_\nu (b_j) - \psi_\nu (b_j)) \|  \\
\| \varphi_\nu (p_\nu) - \varphi_\nu(p_\nu) h \|  \\
\| \varphi_\nu (p_\nu \psi_\nu(h) - p_\nu) \|
\end{array}
\right\} < \frac{\eta^2}{2}\, .  
\end{equation}
Because $\varphi_\nu$ is $n$-decomposable, there are pairwise orthogonal central projections $q^{(0)}, \ldots, q^{(n)} \in F$ (from now on we supress the index $\nu$) satisfying $\be_F  = \sum_l q^{(l)}$ and such that $\ord (\varphi|_{q^{(l)}F}) = 0 \; \forall \, l$. Define $p^{(l)} := q^{(l)} p$. Using ($\ast \ast$), one checks that 
\[
\| \varphi (p^{(l)}) - \varphi (p^{(l)}) h \| < \frac{\eta}{\sqrt{2}}
\]
and, similarly,
\[
\| \varphi(p^{(l)}) - \varphi(p^{(l)}) \varphi(\be_F) \| < \left( \frac{\eta^2}{2} + \frac{\eta^2}{2} \right)^\halb = \eta \, .
\]
Finally, from 
\[
\| \varphi (p \psi(b_j) p - \psi(b_j)) \| < \eta^2 
\]
we see that
\[
\| \varphi( p^{(l)} \psi (b_j) p^{(l)} - q^{(l)} \psi(b_j) q^{(l)}) \| = \| \varphi( q^{(l)} (p \psi(b_j) p - \psi(b_j))) \| < \eta^2 \; \forall \, l  .
\]
\end{nproof}
\alten

\altbn{\label{hereditary}}
\begin{nprop}
Let $A$ be a separable $C^*$-algebra and $B \subset_\her A$ a hereditary $C^*$-subalgebra. Then $\dr B \le \dr A$. 
\end{nprop}

\begin{nproof}
We assume $n := \dr A$ to be finite, for otherwise there is nothing to show. Let $b_1, \ldots , b_m \in B_+$ be normalized elements and $0 < \varepsilon \le 1$; we have to find an $n$-decomposable c.p.\ approximation (for $B$) of $b_1, \ldots, b_m$ within $\varepsilon$. By Remark \ref{r-dr}(ii) we may assume that there is a normalized $h \in B_+$ such that $h b_j = b_j \; \forall \, j$. Set $\eta := \frac{\varepsilon^2}{(24(n+1))^2}$ and choose a c.p.\ approximation $(F, \psi, \varphi)$ with projections $q^{(l)}, \, p^{(l)}$, $l= 0, \ldots, n$ in $F$ as in Lemma \ref{dominated-approximations}.\\
It is then clear that $\ord (\varphi|_{p^{(l)}Fp^{(l)}}) = 0 \; \forall \, l$, so we may apply Lemma \ref{hereditary-perturbation} to obtain c.p.c.\ maps 
\[
\hat{\varphi}^{(l)} : p^{(l)} F p^{(l)} \to \overline{hAh} \subset B
\]
with $\ord \hat{\varphi}^{(l)} = 0$ and $\| \hat{\varphi}^{(l)} (x) - \varphi(x)\| < 8 \eta^\halb =  \varepsilon / (3 (n+1))$ for all $l$ and $x \in (p^{(l)} F p^{(l)})_+$ with $\|x\| \le 1$.\\
Define $\hat{\varphi} : p F p \to B$ by setting
\[
\hat{\varphi} :=  \sum_{l=0}^n \hat{\varphi}^{(l)} \, ,
\]
then $\hat{\varphi}$ is c.p.\ and $n$-decomposable by construction. By \ref{dominated-approximations}(1) the $q^{(l)} \in F$ are central and sum up to $\be_F$ and by \ref{dominated-approximations}(3) we have 
\[
\| \varphi ( q^{(l)} \psi(b_j) q^{(l)} - p^{(l)} \psi(b_j) p^{(l)})\| < \eta \; \forall \, j,l \, .
\]
As a consequence, for each $j$ we obtain
\[
\begin{array}{rcl}
\| \hat{\varphi} \psi (b_j) - b_j \| & \le & \| \hat{\varphi} (\sum_l q^{(l)} \psi (b_j) q^{(l)}) - \hat{\varphi} (\sum_l p^{(l)} \psi (b_j) p^{(l)}) \| \\
& & + \|  \sum \hat{\varphi}^{(l)}  (p^{(l)} \psi (b_j) p^{(l)}) -  \sum \varphi (  p^{(l)} \psi (b_j) p^{(l)}) \| \\
& & + \|  \sum \varphi (  p^{(l)} \psi (b_j) p^{(l)})  -  \sum \varphi (  q^{(l)} \psi (b_j) q^{(l)}) \| \\
& & + \|  \varphi \psi (b_j) - b_j \| \\
& < & (n+1) \cdot \eta + (n+1) \cdot 8 \eta^\halb + (n+1) \cdot \eta + \eta \\
& < & \frac{2}{3} \varepsilon \, .
\end{array}
\]
$\hat{\varphi}$ might fail to be contractive, but by the choice of $\hat{\varphi}^{(l)}$ we have 
\[
\| \hat{\varphi}^{(l)} (p^{(l)}) - \varphi( p^{(l)}) \| < 8 \eta^\halb \, ,
\]
and therefore 
\[
\| \hat{\varphi} \| \le \| {\textstyle \sum} \hat{\varphi}^{(l)}(p^{(l)}) \| \le \|{\textstyle \sum} \varphi (p^{(l)}) \| + (n+1) \cdot 8 \eta^\halb \le 1 + \frac{\varepsilon}{3} \, .
\]
Hence we may replace $\hat{\varphi}$ by the contraction $\frac{1}{1+\frac{\varepsilon}{3}} \cdot \hat{\varphi}$ without disturbing it by more than $\frac{\varepsilon}{3}$. So $(pFp, p \psi(\, . \,)p, \frac{1}{1+\frac{\varepsilon}{3}} \hat{\varphi})$ is a c.p.\ approximation of $B$ for $b_1, \ldots, b_m$ within $\varepsilon$; this completes the proof.
\end{nproof}
\alten

\altbn{\label{matrices-dr=}}
\begin{ncor}
(i) For any separable $C^*$-algebra $A$ and $r \in \N$ we have $\dr A = \dr (M_r \otimes A) = \dr (\Kh \otimes A)$.\\
(ii) If $B \subset_\her A$ is a full hereditary $C^*$-subalgebra, then $\dr B = \dr A$.
\end{ncor}

\begin{nproof}
(i) We have $\dr A \le \dr(M_r \otimes A) \le \dr(\Kh \otimes A)$ by Proposition \ref{hereditary} and $\dr(\Kh \otimes A) \le \dr A$ by \ref{dr-permanence}.\\
(ii) If $B \subset_\her A$ is full, then $B \otimes \Kh \cong A \otimes \Kh$ by Brown's theorem (\cite{Br}, Corollary 2.6) and the assertion follows from (i).
\end{nproof}
\alten

\altbn{\label{continuous-trace}}
\begin{ncor}
Let $A$ be a separable continuous trace $C^*$-algebra. Then $\dr A = \dim \hat{A}$.
\end{ncor}

\begin{nproof}
To show that $\dr A \le \dim \hat{A}$, we modify the proof of \cite{Wi2}, Theorem 4.2.1: There it was shown that $\hr A \le \dim \hat{A}$; in the proof an open covering $(V_1, \ldots, V_l)$ of $\hat{A}$ of strict order $n$ was used to produce a c.p.\ approximation (for $A$) of strict order $n$. Using Proposition \ref{dim=decomposable} one can assume the open covering to be $n$-decomposable; then, replacing the words "of strict order not exceeding $n$" by "$n$-decomposable", literally the same construction as in \cite{Wi2} yields an $n$-decomposable c.p.\ approximation for $A$, so $\dr A \le \dim \hat{A}$.\\
Conversely, note that $A$ is locally stably isomorphic to $\Ch_0(\hat{A})$, more precisely: each $t \in \hat{A}$ has a compact neighborhood $M_t$ such that $A/J_{M_t} \otimes \Kh \cong \Ch(M_t) \otimes \Kh$. Using Proposition \ref{dim=dr} and Corollary \ref{matrices-dr=}, we have
\[
\dim M_t = \dr \Ch(M_t) = \dr (\Ch(M_t) \otimes \Kh) = \dr (A/J_{M_t} \otimes \Kh) \le \dr A \, .
\]
Since $A$ is separable, $\hat{A}$ is covered by countably many sets $M_t$ with $\dim M_t \le \dr A$. Now the countable sum theorem (cf.\ \cite{HW}, Theorem III.2) for covering dimension says that $\dim \hat{A} \le \dr A$. Note that the estimate does not follow from $\hr A \le \dr A$, since for the homogeneous rank we only have a slightly weaker result. 
\end{nproof}
\alten

\altbn{\label{unitization}}
\begin{nprop}
Let $A$ be a separable $C^*$-algebra and $\tilde{A}$ its smallest unitization. Then $\dr A = \dr \tilde{A}$.
\end{nprop}

\begin{nproof}
$A \lhd \tilde{A}$ is an ideal, so $\dr A \le \dr \tilde{A}$ by Proposition \ref{hereditary}; we show $\dr \tilde{A} \le \dr A =: n$ (again we only have to prove the statement for $n < \infty$).\\
Suppose $\varepsilon > 0$ and $b_1, \ldots, b_m \in A_+$ with $\|b_j\| \le 1$ are given. We have to find an $n$-decomposable c.p.\ approximation of $\tilde{A}$ for $\be, b_1, \ldots, b_m$ within $\varepsilon$, where $\be$ denotes the unit of $\tilde{A}$. As in Remark \ref{r-dr}(ii), we may assume that there is $h \in A_+$ with $\|h\| \le 1$ and $h b_j = b_j \; \forall \, j$. \\
Choose $\gamma > 0$ such that $\gamma + 12 \gamma^\halb < \varepsilon$. Since orthogonality is a weakly stable relation (cf.\ \ref{weakly-stable} or \cite{Lo}, 10.1.10), there is $\eta > 0$ such that, whenever $e_0, \, e_1$ are positive normalized elements in $\tilde{A}$ with $\|e_0 e_1\| < \eta$, there are positive normalized $d_0, \, d_1 \in \tilde{A}$ with $d_0 \perp d_1$ and $\|e_k - d_k\| < \gamma$, $k = 0, 1$. We may assume $\eta$ to be smaller than $\halb (\varepsilon - 12 \gamma^\halb)$. Now apply Lemma \ref{dominated-approximations} to obtain a c.p.\ approximation $(F,\psi, \varphi)$ (of $A$) for $h, b_1, \ldots, b_m$ within $\eta$ and projections $p^{(l)} \le q^{(l)} \in F$, $l=0, \ldots,n$, such that (1), (2) and (3) of \ref{dominated-approximations} hold (in particular the $q^{(l)}$ are pairwise orthogonal, central, and sum up to $\be_F$). Set 
\[
q := p^{(0)} + {\textstyle \sum_1^n} q^{(l)}
\]
and
\[
\tilde{F} := q F q \oplus \C = p^{(0)} F p^{(0)} \oplus ({\textstyle \bigoplus_{l=1}^n } q^{(l)} F) \oplus \C
\]
and define a u.c.p.\ map $\tilde{\psi} : \tilde{A} \to \tilde{F}$ as the unitization of $\psi_q : A \to q F q$ (using \cite{CE}, Lemma 3.9); then 
\[
\tilde{\psi}(a) = q \psi(a) q \; \forall \, a \in A \, \mbox{ and } \, \tilde{\psi}(\be) = q \oplus \be_\C \, .
\]
Now (using that $\varphi(p^{(0)}) \perp \varphi(q^{(0)} - p^{(0)})$)
\begin{eqnarray*}
\| \varphi(p^{(0)}) (\be - \varphi(q))\| & = & \| \varphi(p^{(0)}) (\be - \varphi(\be_F) + \varphi(q^{(0)} - p^{(0)})) \| \\
& = & \| \varphi(p^{(0)}) (\be - \varphi(\be_F)) \| \\
& \stackrel{\ref{dominated-approximations} (3)}{<} & \eta \, ,
\end{eqnarray*}
so by the choice of $\eta$ there are $0 \le d_0, d_1 \in \tilde{A}$ with $\| d_k \| \le 1$, $d_0 \perp d_1$ and $\| \varphi(p^{(0)}) - d_0 \|$, $\|(\be - \varphi(q)) - d_1 \| < \gamma$.\\
Next observe that 
\[
\| \varphi(p^{(0)}) (\be - g_{0,\gamma}(d_0)) \| \le \| d_0 (\be - g_{0, \gamma}(d_0)) \| + \gamma \le 2 \gamma \, 
\]
(cf.\ \ref{functions} for the definition of $g_{\alpha,\beta}$), so we may apply Lemma \ref{hereditary-perturbation} to obtain
\[
\hat{\varphi} : p^{(0)} F p^{(0)} \to \overline{g_{0,\gamma}(d_0) A g_{0,\gamma}(d_0)} = \overline{d_0 A d_0}
\]
with $\ord \hat{\varphi} = 0$ and 
\[
\| \hat{\varphi}(x) - \varphi(x) \| < 8 \cdot (2 \gamma)^\halb < 12 \gamma^\halb
\]
for all $x \in (p^{(0)} F p^{(0)})_+$ with $\|x\| \le 1$. Define
\[
\tilde{\varphi} : \tilde{F} \to \tilde{A}
\]
by
\[
\tilde{\varphi}|_{p^{(0)} F p^{(0)}} := \hat{\varphi}\, , \; \tilde{\varphi}|_{\bigoplus_{l=1}^n q^{(l)} F} := \varphi|_{\bigoplus_{l=1}^n q^{(l)} F} \, \mbox{ and } \; \tilde{\varphi}(\be_\C) := d_1 \, .
\]
By construction, $\tilde{\varphi}|_{p^{(0)} F p^{(0)} \oplus \C}$ and $\tilde{\varphi}|_{q^{(l)}F}$, $l=1, \ldots, n$, all have strict order zero, so $\tilde{\varphi}$ is $n$-decomposable. Furthermore, we have
\begin{eqnarray*}
\| \tilde{\varphi} \tilde{\psi} (\be) - \be \| & = & \| \hat{\varphi}(p^{(0)}) + \varphi( {\textstyle \sum_{l=1}^n } q^{(l)}) + d_1 - \be \| \\
& \le & \| \varphi(q) + d_1 - \be \| + 12 \gamma^\halb \\
& \le & \gamma + 12 \gamma^\halb < \varepsilon
\end{eqnarray*}
and, because $\| \varphi( p^{(0)} \psi(b_j) p^{(0)} - q^{(0)} \psi(b_j) q^{(0)}) \| < \eta $ by \ref{dominated-approximations}(3),
\begin{eqnarray*}
\| \tilde{\varphi} \tilde{\psi} (b_j) - b_j\| & = & \| \hat{\varphi}(p^{(0)} \psi(b_j) p^{(0)}) + \varphi( {\textstyle \sum_{l=1}^n } q^{(l)} \psi(b_j)) - b_j \| \\
& \le & \| \varphi \psi (b_j) - b_j \| + 12 \gamma^\halb + \eta \\
& \le & 12 \gamma^\halb + 2 \eta < \varepsilon \, .
\end{eqnarray*}
Again it might happen that $\tilde{\varphi}$ is not contractive, but since $| 1 - \|\tilde{\varphi} \| | < \varepsilon$, multiplication by $\frac{1}{1+ \varepsilon}$ yields a contraction which differs from $\tilde{\varphi}$ by at most $\varepsilon$.
\end{nproof}
\alten

\section{Quasidiagonality}

\altbn{\label{qd}}
Let $A$ be a $C^*$-algebra and $\pi: A \to \Bh(\Hh)$ a representation. Recall from \cite{Vo2} that $\pi$ is said to be quasidiagonal if there is a sequence $(p_n)_\N$ of finite-rank projections in $\Bh(\Hh)$ converging strongly to $\be_\Hh$ and satisfying $\| [p_n,\pi(a)]\| \to 0 \; \forall \; a \in A$. $A$ is quasidiagonal as a $C^*$-algebra if it has a faithful quasidiagonal representation; $A$ is strongly quasidiagonal if every representation is quasidiagonal.\\
We will proceed to show that, if $\dr A$ is finite, $A$ is quasidiagonal; since decomposition rank passes to quotients, it will even follow that $A$ is strongly quasidiagonal (which in particular implies that $A$ is inner quasidiagonal in the sense of \cite{BK2}).
\alten

\altbn{\label{approximately-multiplicative}}
\begin{nprop}
A separable $C^*$-algebra $A$ has decomposition rank less than or equal to $n$, if and only if the following holds:\\
For any $b_1, \ldots, b_m \in A$ and $\varepsilon > 0$ there is a c.p.\ approximation $(F,\psi, \varphi)$ such that
\[
\| \varphi \psi (b_k) - b_k \| , \, \| \psi(b_k) \psi(b_l) - \psi (b_k b_l) \| < \varepsilon \; \forall \, k, l
\]
and such that $\varphi$ is $n$-decomposable.\\
In particular, if $\dr A < \infty$, then $A$ is quasidiagonal.
\end{nprop}

\begin{nproof}
We may assume $0 \le b_k \le \be \; \forall \, k$ and $\varepsilon \le 1$; by Lemma \ref{multiplicative-domain} it suffices to find a c.p.\ approximation $(F,\psi, \varphi)$ for $b_1, \ldots, b_m$ within $\varepsilon$ such that $\| \psi(b_k^2) - \psi(b_k)^2\| < \varepsilon^2$.\\
Let $(F', \psi', \varphi')$, $F' = M_{r_1} \oplus \ldots \oplus M_{r_s}$, be a c.p.\ approximation for $b_1, \ldots, b_m$, $b_1^2, \ldots,b_m^2$ within $\delta := \frac{\varepsilon^4}{6 (n+1)}$ such that $\varphi'$ is $n$-decomposable. We write $\psi_i'$, $\varphi_i'$ and $\be_i$ for the respective summands of $\psi'$, $\varphi'$ and $\be_F$. Note that $\ord \varphi_i' = 0 \; \forall \, i$ and that (by Proposition \ref{strict-order} and since $C^*$-norms are cross-norms),
\[
\|\varphi_i'(x)\| = \| \pi_\varphi( \id_{(0,1]} \otimes x) \| = \|x\| \cdot \|\varphi_i'(\be_i)\| \; \forall \; x \in M_{r_i} \, .
\]
Define 
\[
I:= \{ i \in \{1, \ldots,s\} \, | \, \|\psi_i(b_k^2) - \psi_i(b_k)^2 \| \ge \varepsilon^2 \mbox{ for some } k \} \, .
\]
Thus, for each $j \in I$ we have (for a suitable $k$)
\begin{eqnarray*}
\varepsilon^2 \cdot \| \varphi_{j}'(\be_{j})\| & \le & \| \varphi_{j}' (\psi_{j}' (b_k^2) - \psi_{j}'(b_k)^2) \| \\
& \le & \| {\textstyle \sum_i} \varphi_i' (\psi_i' (b_k^2) - \psi_i'(b_k)^2) \| \\
& = & \| \varphi' (\psi' (b_k^2) - \psi'(b_k)^2) \| \\
& \le & \| \varphi' \psi' (b_k^2) - \varphi' \psi' (b_k)^2 \| \\
& \le & \| \varphi' \psi' (b_k^2) - b_k^2\| + \| b_k^2 - \varphi' \psi' (b_k)^2 \| \\
& < & 3 \delta \, ,
\end{eqnarray*}
so $\| \varphi'(\be_{j})\| < \frac{3 \delta}{\varepsilon^2}$ for all ${j} \in I$. But now by Lemma \ref{l-decomposable}, we have $\| \sum_{i \in I} \varphi'(\be_i) \| \le (n+1) \cdot \frac{3 \delta}{\varepsilon^2}$. Define $(F, \psi, \varphi)$ by setting $F:= \bigoplus_{i \in \{1, \ldots,s\} \setminus I} M_{r_i} \subset F'$, $\psi(a):= \be_F \psi'(a) \be_F$ and $\varphi:= \varphi'|_F$, then 
\[
\| \varphi \psi (b_k) - b_k\|, \, \|\varphi \psi (b_k^2) - b_k^2 \| < \delta + \frac{(n+1) 3 \delta}{\varepsilon^2} = \frac{\varepsilon^4}{6(n+1)} + \frac{\varepsilon^2}{2} < \varepsilon^2 \, .
\]
The last statement follows from \cite{BK1}, Theorem 5.2.2.
\end{nproof}
\alten

\altbn{\label{dr-nf}}
\begin{nremarks}
(i) We do not have a similar statement for the completely positive or the homogeneous rank, but at least one can show that, if $A$ is simple, unital and if $\cpr A \le 1$, then $A$ is quasidiagonal.\\
(ii) If $\dr A = n$ and $p_1, \ldots, p_m$ are projections in $A$, we can choose our c.p.\ approximation $(F, \psi, \varphi)$ in a way such that the $\psi(p_k)$ are close to projections in $F$. In particular, if $A$ is unital we may assume $\psi(\be_A)$ to be close to a projection $q \in F$. But then the unital map $\psi': A \to qFq$ defined by 
\[
a \mapsto (q \psi(\be_A) q)^{-\halb} \psi(a) (q \psi(\be_A) q)^{-\halb}
\]
is close to $\psi$. Thus, setting $F' := qFq$ and $\varphi' := \varphi|_{qFq}$, we have a c.p.\ approximation $(F, \psi', \varphi')$ such that $\varphi'$ is $n$-decomposable (and close to being unital) and $\psi'$ is unital. (The ``unitized'' c.p.\ approximation $(F \oplus \C, \tilde{\psi }, \tilde{\varphi })$ will in general only be $(n+1)$-decomposable, but cf.\ Proposition \ref{unitization}.)
\end{nremarks}
\alten

\altbn{\label{t-qd}}
\begin{ntheorem}
Let $A$ be a separable $C^*$-algebra with $\dr A = n < \infty$. Then $A$ is strongly quasidiagonal.
\end{ntheorem}

\begin{nproof}
By \cite{BK2}, Proposition 2.8, it suffices to show that any irreducible representation $\pi : A \to \Bh(\Hh)$ is quasidiagonal. But $\dr \pi(A) \le \dr A$ by \ref{dr-permanence}(iii), so $\pi(A)$ has a faithful quasidiagonal representation $\rho : \pi(A) \to \Bh(\Hh')$ by Proposition \ref{approximately-multiplicative}. Now if $\pi(A)$ doesn't meet the compacts, then the identity representation of $\pi(A)$ is approximately unitarily equivalent to $\rho$ by Voiculescu's theorem (\cite{Vo1}), hence $\pi$ is quasidiagonal.

Next suppose $\pi(A) \cap \Kh(\Hh) \neq \{0\}$, then $\Kh(\Hh) \subset \pi(A)$ since $\pi$ is irreducible. We will show the following:\\
For an arbitrarily small $\varepsilon > 0$, a rank-one projection $a_0 \in \Kh(\Hh)$ and $a_1, \dots, a_m \in \pi(A)_+$ with $\|a_j\| \le 1$ there is a projection $p \in \Kh(\Hh)$ such that $\|p a_0 p \| > 1 - \varepsilon$ and $\|p a_j - a_j p\| < \varepsilon \; \forall \; j$.\\
It then follows that $\pi$ is indeed quasidiagonal, for if $\{a_1, \ldots, a_m\}$ contains enough partial isometries with the same source projection $a_0$, the condition $\|p a_j - a_j p\| < \varepsilon$ will guarantee that $p$ can be chosen arbitrarily close to $\be_\Hh$ in the strong operator topology.

So choose $\delta>0$ with $4 (n+1) \delta^\halb + 2 \delta < \varepsilon$ and a c.p.\ approximation $(F= M_{r_1} \oplus \ldots \oplus M_{r_s}, \psi, \varphi)$ (of $\pi(A)$) for $a_0, a_1, \ldots, a_m, a_1^2, \ldots, a_m^2$ within $\delta$ such that $\varphi$ is $n$-decomposable. By Remark \ref{dr-nf}(ii) we may assume $\psi(a_0)$ to be close to a projection; in particular we can assume $\| \be_i \psi(a_0)\|$ to be smaller than $\delta$ or larger than $1-\delta$ for all $i$. Now from Lemma \ref{l-decomposable} we see that there is some $i \in \{1, \ldots, s\}$ (which we fix from now on) with $\|\be_i \psi(a_0)\| > 1-\delta$ and $\|\varphi(\be_i \psi(a_0))\| > \frac{1-\delta}{n+1}$.\\
Set $C:= C^*(\varphi(M_{r_i})) \subset \pi(A)$ and let $\rho: C \to M_{r_i}$ be the unique representation with $\rho \verk \varphi (x) = \|\varphi(\be_i)\| \cdot x \; \forall \, x \in M_{r_i}$ (by Proposition \ref{strict-order}, $C$ is a quotient of $\Ch_0((0,1]) \otimes M_{r_i}$). $\rho$ is irreducible and can be extended to all of $A$, i.e., there are an irreducible representation $\bar{\rho} : \pi(A) \to \Bh(\bar{\Hh})$ and an isometry $V: \C^{r_i} \to \bar{\Hh}$ with $\rho(x) = V^* \bar{\rho}(x) V \; \forall \, x \in C$. In particular we have $\bar{\rho}(\Kh(\Hh)) \neq \{0\}$, so by \cite{KR}, Theorem 10.4.6, there is a unitary $U : \bar{\Hh} \to \Hh$ such that $U \bar{\rho}(a) U^* = a \; \forall \; a \in \pi(A)$. Define $p := UVV^*U^* \in \Kh(\Hh)$. Using that $\rho$ is a $*$-homomorphism, it is straightforward to check that $[p,C] = \{0\}$ and that $p \varphi(\be_i) = \|\varphi(\be_i)\| \cdot p$.\\
Now by Lemma \ref{multiplicative-domain} and because $\varphi \verk \psi$ approximates $a_0$ within $\delta$ we obtain 
\begin{eqnarray*}
\| \varphi(\be_i \psi(a_0)) - \varphi(\be_i) \varphi \psi(a_0)\| & \stackrel{\ref{multiplicative-domain}}{\le} & \| \varphi(\be_i)\| \| \varphi(\psi(a_0)^2) - \varphi \psi(a_0)^2 \|^\halb \\
& \le & \| \varphi(\be_i)\| \| \varphi \psi(a_0^2) - \varphi \psi(a_0)^2 \|^\halb \\
& \le & \| \varphi(\be_i)\| (3 \delta)^\halb \, .
\end{eqnarray*}
We now have
\begin{eqnarray*}
\|pa_0p\| & \ge & \|p \varphi \psi(a_0) p\| - \delta \\
& = & \frac{1}{\|\varphi(\be_i)\|} \|p \varphi(\be_i) \varphi \psi(a_0) p\| - \delta\\
& \ge & \frac{1}{\|\varphi(\be_i)\|} \|p \varphi(\be_i \psi(a_0)) p\| - \frac{(n+1) \cdot 3^\halb \delta^\halb}{1-\delta} - \delta \\
& = & \frac{1}{\|\varphi(\be_i)\|} \|\rho \verk \varphi(\be_i \psi(a_0)) \| - \frac{(n+1) \cdot 3^\halb \delta^\halb}{1-\delta} - \delta \\
& = & \| \be_i \psi(a_0)\| - \frac{(n+1) \cdot 3^\halb \delta^\halb}{1-\delta} - \delta \\
& > & 1- 2 \delta - \frac{(n+1) \cdot 3^\halb \delta^\halb}{1-\delta} \\
& > & 1 - \varepsilon \, .
\end{eqnarray*}
Similarly,
\begin{eqnarray*}
\| p a_j - a_j p \| & \le & \frac{1}{\|\varphi(\be_i)\|} \| p \varphi(\be_i) \varphi \psi (a_j) - \varphi \psi (a_j) \varphi(\be_i) p \| + 2 \delta \\
& \le & \frac{1}{\|\varphi(\be_i)\|} \| p \varphi(\be_i \psi (a_j) - \psi (a_j) \be_i) p \| + \frac{2(n+1) \cdot 3^\halb \delta^\halb}{1-\delta} + 2 \delta \\
& = & 0 + \frac{2(n+1) \cdot 3^\halb \delta^\halb}{1-\delta} + 2 \delta \\
& < & \varepsilon
\end{eqnarray*}
for $j= 0, \ldots, m$ and we are done.
\end{nproof}
\alten

\section{Quasidiagonal extensions}

Given an extension $0 \to J \to A \to A/J \to 0$ of $C^*$-algebras, we already know from \ref{dr-permanence} that $\max \{ \dr J, \dr A/J \} \le \dr A$. The Toeplitz extension $0 \to \Kh \to \Th \to \Ch(S^1) \to 0$ shows that one does not have equality in general ($\dr \Th = \infty$, since $\Th$ is not quasidiagonal). It is therefore natural to ask for conditions under which $\dr A$ is determined by $\dr J$ and $\dr A/J$. A class of extensions which behave particularly well in this respect are the quasidiagonal ones; these were introduced in \cite{Sal}.

\altbn{\label{qd-extensions}}
\begin{nprop}
Let $A$ be a $C^*$-algebra and $J \lhd A$ an ideal  with a quasicentral approximate unit consisting of projections. Then $\dr A = \max \{ \dr J, \dr A/J \}$.
\end{nprop}

\begin{nproof}
Since $n:= \max \{ \dr J, \dr A/J \} \le \dr A$, we only have to show $\dr A \le n$ for $n < \infty$.\\
Given $b_1, \ldots, b_m \in A_+$ with $\|b_i\| \le 1$ and $\varepsilon > 0$, choose a c.p.\ approximation $(F,\psi, \varphi)$ (of $A/J$) for $\pi(b_1), \ldots, \pi(b_m)$ within $\frac{\varepsilon}{4}$ such that $\varphi$ is $n$-decomposable; here, $\pi : A \to A/J$ denotes the quotient map. By \ref{liftable}, one can lift $\varphi$ to an $n$-decomposable c.p.\ (but not necessarily contractive) map $\hat{\varphi} : F \to A$. Choose $\delta > 0$ such that the assertion in \ref{weakly-stable} holds with $\eta := \frac{\varepsilon}{4}$. Take a set $\{e_1, \ldots, e_k\}$ of matrix units for $F$. Using our assumption it is routine to find a projection $p \in J$ such that the following hold:\\
(i) The map $\varphi' : F \to (\be-p)A(\be-p)$, given by $\varphi'(\, . \,) := (\be-p) \hat{\varphi} (\, .\,) (\be-p)$, satisfies  $\| \varphi'\| < 1 + \delta$.\\ 
(ii) For $j = 1, \ldots ,k$, the elements $\varphi'(e_j) \in (\be-p)A(\be-p)$ satisfy the relations $(\Rh)$ defining $n$-decomposability up to $\delta$ (cf.\ \ref{weakly-stable}).\\
(iii) $\| p b_i p + (\be - p) b_i (\be - p) - b_i\| < \frac{\varepsilon}{4}$.\\
(iv) $\| \varphi'(\psi(b_i)) - (\be - p) b_i (\be - p)\| < \frac{\varepsilon}{4}$. \\
But now by \ref{weakly-stable} there is an $n$-decomposable c.p.c.\ map $\varphi'' : F \to (\be -p) A (\be-p)$ with $\| \varphi'' - \varphi' \| < \frac{\varepsilon}{4}$.\\
Since $pAp \subset_\her J$, we have $\dr pAp \le n$ by Proposition \ref{hereditary} and may choose a c.p.\ approximation $(\bar{F}, \bar{\psi}, \bar{\varphi})$ (of $pAp$) for $p b_i p$, $i =1, \ldots, m$ within $\frac{\varepsilon}{4}$ and with $\bar{\varphi}$ $n$-decomposable.\\
Now $\varphi'' \oplus \bar{\varphi} : F \oplus \bar{F} \to (\be-p)A(\be-p) \oplus pAp \subset A$ is c.p.c.\ and $n$-decomposable, since $\varphi''$ and $\bar{\varphi}$ are both c.p.c.\ and $n$-decomposable and have orthogonal images.
The map $\bar{\psi}: pAp \to \bar{F}$ has a c.p.c.\ extension to all of $A$, defined by $a \mapsto \bar{\psi}(pap)$; this extension we also denote by $\bar{\psi}$. It only remains to check that $(F \oplus \bar{F}, \psi \oplus \bar{\psi}, \varphi'' \oplus \bar{\varphi})$ approximates $b_1, \ldots, b_m$ within $\varepsilon$. But for each $i$ we have
\begin{eqnarray*}
\| (\varphi'' \oplus \bar{\varphi}) \verk (\psi \oplus \bar{\psi}) (b_i) - b_i\| & = & \| \varphi'' \psi (b_i) + \bar{\varphi} \bar{\psi} (p b_i p) - b_i \| \\
& \le & \| \varphi'' \psi (b_i) - \varphi' \psi (b_i) \| \\
& & + \| \varphi' \psi(b_i) - (\be-p) b_i (\be-p) \| \\
& & + \|\bar{\varphi} \bar{\psi} (p b_i p) - p b_i p \| \\
& & + \| (\be-p) b_i (\be-p) + p b_i p - b_i \| \\
& < & 4 \cdot \frac{\varepsilon}{4} = \varepsilon
\end{eqnarray*}
and the proof is complete.
\end{nproof}
\alten

\altbn
If $J = \Kh$, $\Kh$ being the compact operators on a separable Hilbert space, then quasidiagonal extensions are precisely those which behave well with respect to the decomposition rank:

\begin{nprop}
Let $A$ be a separable $C^*$-algebra containing $\Kh$ as an ideal and with $\dr A/\Kh < \infty$. Then $\dr A = \dr A/\Kh$ if $\Kh$ contains a quasicentral approximate unit consisting of projections; otherwise, $\dr A = \infty$.
\end{nprop}

\begin{nproof}
Suppose $\dr A$ and $\dr A/\Kh$ are both finite. We only have to show that $\Kh$ contains a quasicentral approximate unit consisting of projections, then the assertion will follow from Proposition \ref{qd-extensions}.\\
But $A$ is strongly quasidiagonal by Theorem \ref{t-qd}, so there is a quasidiagonal irreducible representation $\pi : A \to \Bh(\Hh)$ with $\pi(\Kh) \neq 0$. $\pi$ restricts to a nonzero irreducible representation of $\Kh$ on $\pi(\Kh) \Hh \subset \Hh$. Now since $\Kh$ has, up to unitary equivalence, only one irreducible representation, $\pi(\Kh)$ coincides with the compact operators on $\pi(\Kh) \Hh$, and these are contained in $\Kh(\Hh)$. In particular we see that $\pi(\Kh) \cap \Kh(\Hh) \neq 0$, so $\Kh(\Hh) \subset \pi(A)$ since $\pi$ is irreducible. As a consequence, $\pi(\Kh)$ is a nonzero ideal of $\Kh(\Hh)$, so $\pi(\Kh) = \Kh(\Hh)$. \\
Choose a quasicentral (with respect to $\pi(A)$) approximate unit of projections $(q_k)_{\N}$ for $\Kh(\Hh)$; this is possible, since $\pi(A)$ acts quasidiagonally on $\Hh$. But then $(\pi^{-1}(q_k))_{\N}$ is an approximate unit of projections for $\Kh \subset A$. To see that it is quasicentral for $A$, note that for any $a \in A$ we have $[\pi^{-1}(q_k),a] \in \Kh \; \forall \, k$, hence
\[
\| [\pi^{-1}(q_k),a] \| = \| \pi( [\pi^{-1}(q_k),a] ) \| = \| [q_k, \pi(a)] \| \to 0 \, .
\]
\end{nproof}
\alten

\section{Some examples and conclusive remarks}

\altbn{\label{examples}}
We finally give a list of examples; see \cite{Bl1}, \cite{Da} or \cite{Ro} for detailed information on the construction and properties of these.

\begin{nexamples}
(i) It is clear from \ref{dr-permanence} that, if $A$ is an $AF$ algebra, then $\dr A = 0$. Conversely, if $\dr A = 0$, then $\cpr A = 0$ by \ref{r-dr}(i), hence $A$ is $AF$ by \cite{Wi1}, Theorem 4.2.3.\\
(ii) If $A$ is an approximately homogeneous algebra, then $\dr A$ is bounded by the dimensions of the base spaces; this also follows from \ref{dr-permanence}. In particular, for a limit circle algebra $A$ we have $\dr A \le 1$; for the irrational rotation algebras $A_\theta$ we obtain $\dr A_\theta =1$, since these are $AT$, but not $AF$. Furthermore, it is not hard to see from the construction that Blackadar's simple unital projectionless algebra has decomposition rank one.\\
(iii) In \cite{Wi3} it will be shown that, for a subhomogeneous algebra $A$, we have $\dr A = \max_k \dim (\Prim_k A)$, where $\Prim_k A$ is the (locally compact Hausdorff) space of kernels of $k$-dimensional irreducible representations. In particular, the dimension drop intervals used as building blocks for approximately subhomogeneous $C^*$-algebras have decomposition rank one.\\
(iv) Using recent work of Phillips and Q.\ Lin, it follows from (iii) that $\dr (\Ch(M) \rtimes_\alpha \Z) \le \dim M$ if $M$ is a compact smooth manifold and $\alpha$ is a minimal diffeomorphism of $M$.\\
(v) The Toeplitz algebra $\Th$ and the Cuntz algebras $\Oh_n$ are not quasidiagonal, since they contain infinite projections; it  therefore follows from Proposition \ref{approximately-multiplicative} that $\dr \Th = \dr \Oh_n = \infty$.\\
(vi) Recently, R{\o}rdam has given an example of a (non-simple) $C^*$-algebra which is an inductive limit of algebras of the form $\Ch_0 ([0,1)) \otimes M_r$ and absorbs the Cuntz algebra $\Oh_\infty$ tensorially (see \cite{Ro2}). This example is purely infinite in the sense of \cite{KiR} and has decomposition rank one. 
\end{nexamples}
\alten

\altbn
There are several other generalizations of covering dimension to $C^*$-algebras, among which the stable rank (introduced by Rieffel in \cite{Ri1}) and the real rank (introduced by Brown and Pedersen in \cite{BP}) are of particular interest. Both notions behave rather different from the decomposition rank. In particular, they take their lowest values not only on $AF$ algebras (e.g., simple purely infinite $C^*$-algebras always have real rank zero); furthermore, they break down under taking matrix algebras, whereas the decomposition rank is strongly Morita invariant by \ref{matrices-dr=}. Although in all our examples the decomposition rank dominates the real and the stable rank, at the present stage we do not have a general statement which compares the three concepts. 
\alten

\altbn
As we pointed out earlier, one of the reasons for introducing yet another noncommutative version of covering dimension was, to find conditions under which generalized inductive limit $C^*$-algebras are classifiable by $K$-theory data; see \cite{Ro} for an introduction to the classification program.      In subsequent work the following partial result will be proved: If $A$ is separable, simple and unital with real rank zero, stable rank one, finite decomposition rank, weakly unperforated $K_0$-group and a unique tracial state, then $A$ has tracial rank zero in the sense of \cite{Li}. It then follows from Lin's work that $C^*$-algebras which satisfy the preceding conditions and the universal coefficient theorem (cf.\ \cite{Bl1}) are completely classified by their $K$-groups. This result is promising since it shows that, at least under suitable extra conditions, $C^*$-algebras with finite decomposition rank are accessible to classification. 
\alten

\end{document}